# Fourier series (based) multiscale method for computational analysis in science and engineering:

# III. Fourier series multiscale method for linear differential equation with constant coefficients

Weiming Sun[a*+] and Zimao Zhang[b]

**Abstract:** Fourier series multiscale method, a concise and efficient analytical approach for multiscale computation, will be developed out of this series of papers. In the third paper, the analytical analysis of multiscale phenomena inherent in the $2r$-th order linear differential equations with constant coefficients and subjected to general boundary conditions is addressed. The limitation of the algebraical polynomial interpolation based composite Fourier series method is first discussed. This leads to a new formulation of the composite Fourier series method, where homogeneous solutions of the differential equations are adopted as interpolation functions. Consequently, the theoretical framework of the Fourier series multiscale method is provided, in which decomposition structures of solutions of the differential equations are specified and implementation schemes for application are detailed. The Fourier series multiscale method has not only made full use of existing achievements of the Fourier series method, but also given prominence to the fundamental position of structural decomposition of solutions of the differential equations, which results in perfect integration of consistency and flexibility of the Fourier series solution.

**Keywords:** Fourier series, multiscale, differential equation

## 1. Introduction

The researches on the analytical solution of linear differential equation(s) with constant coefficients have been serving as an important motivation for the development of the Fourier series methods for a long time. In particular, the transverse elastic bending problem of thin plates with simply supported boundary conditions, known as the typical fourth order linear differential equation with constant coefficients, gave an example of successful application in the early stage of the Fourier series method [1]. This application not only verified the validity of the theoretical framework of the Fourier series method, but also led to two traditional, specific computational procedures, namely the Navier's method and the Levy's method [2, 3]. In the next two centuries, especially in the last half century, numbers of researchers such as Timoshenko [2], Zhang [3], Gorman [4-22], Huang [23-34], Chaudhuri [1, 35-40], Kabir [41-51], Oktem [52-61], Yan [62] and Li [63-76] successively performed investigations on the boundary value problems of the linear differential equation or highly coupled linear differential equations, which are with constant coefficients and subject to general boundary conditions, and proposed several more general, systematic computational procedures based on the theoretical framework of the Fourier series method. These procedures include the Fourier series based superposition method, the Fourier series based general analytical method, the Fourier series direct-expansion method and the Fourier series method with supplementary terms, etc, which combine to greatly expand the range of application of the Fourier series method and gain incomparable advantages (convenience for computation, simplicity in formulation) over numerical methods. In recent years, some multiscale phenomena, such as the localized high oscillation, sharp gradients or singularity, inherent in the linear differential equation(s) with constant coefficients have drawn close attention from researchers [77-80]. Accordingly, with a goal of uniform approximation and accurate capture of the localized multiscale structures within the solution domain, the development of highly accurate multiscale methods in the theoretical framework of analytical methods has recently become the important issue of the innovation of the Fourier series method.

[a]Department of Mathematics and Big Data, School of Artificial Intelligence, Jianghan University, Wuhan, 430056, China
[b]Department of Mechanics, School of Civil Engineering, Beijing Jiaotong University, Beijing, 100044, China
*Correspondence to: Weiming Sun, Department of Mathematics and Big Data, School of Artificial Intelligence, Jianghan University, Wuhan, 430056, China
[+] E-mail: xuxinenglish@hust.edu.cn

  Fortunately, by taking the composite Fourier series method for simultaneous approximation to functions and their (partial) derivatives that arose in the second paper [81] of the series researches on the Fourier series multiscale method as the basis, together with the usually used Fourier coefficient comparison method, or other derivation methods for discrete equations such as the numerical methods in strong form or the numerical methods in weak form, we have actually constructed a set of systematic method to analytically solve the $2r$-th ($r$ is a positive integer) order linear differential equation(s) with constant coefficients, herein also called after the composite Fourier series method. For the composite Fourier series method, the obtained series solution maintains a consistency of the structure in some sense, which means that the first priority has been given to the satisfaction of the sufficient conditions for the $2r$ times term-by-term differentiation of the Fourier series solution. And consequently, it becomes a routine task for the decomposition of the solution into a corner function that describes the discontinuities at corners of the domain (only for the two-dimensional function), boundary functions that describe the discontinuities at boundaries of the domain and an internal function that describes the smoothness within the domain. Meanwhile, the obtained series solution possesses the flexibility of the formulation to a certain degree and allows of further adjustments in expression. Usually the corner function and the internal function can be selected respectively as the algebraical polynomials with specific orders and Fourier series. However, the selection of interpolation basis functions for the boundary functions might differ in different situations. Specifically, in the process of solving the linear differential equation(s) with constant coefficients, the selected interpolation basis functions satisfy the homogeneous form of the differential equation(s) to be solved, and then the boundary functions constitute the general solutions of the original linear differential equation(s) with constant coefficients. This also means that the existing research findings through centuries of the Levy's method, and especially the Fourier series based general analytical method, can be adopted for the construction of boundary functions. Therefore, the newly developed composite Fourier series method that is based on the interpolation with homogeneous solutions of the differential equation(s) has not only made full use of existing achievements of the Fourier series method, but also given prominence to the fundamental position of structural decomposition of solutions of the differential equations, which results in perfect integration of consistency and flexibility of the Fourier series solution and can be applied directly to the analytical analysis of multiscale phenomena inherent in the linear differential equation(s) with constant coefficients. This multiscale method stems from the general framework of the Fourier series method and thus it bears the name of Fourier series multiscale method.

For brevity, in the third paper of the series of researches on Fourier series multiscale method, we confine ourselves to a systematic development of the Fourier series multiscale method of the linear differential equation with constant coefficients. However, the standpoints and techniques used herein are applicable to the highly coupled linear differential equations with constant coefficients and subjected to general boundary conditions. In this paper, the limitation of the algebraical polynomial interpolation based composite Fourier series method is first demonstrated. Consequently, the implementation schemes for application of the Fourier series multiscale method, such as structural decomposition of the solution, derivation of the general solution, equivalent transformation of the solution, derivation of the supplementary solution, and derivation of discrete equations, are detailed. Then the multiscale characteristic of the Fourier series multiscale solution is briefly discussed. Numerical examples of certain classes of linear differential equation(s) with constant coefficients, such as the convection-diffusion-reaction equation, elastic bending of a thick plate on biparametric foundation, and wave propagation in an infinite rectangular beam, will be presented successively in the forthcoming papers of this series of researches on the Fourier series multiscale method.

## 2. Limitation of algebraical polynomial interpolation based composite Fourier series method

In [81], we have validated the algebraical polynomial interpolation based composite Fourier series method with four different kinds of sample functions, i.e., one-dimensional polynomial, one-dimensional trigonometrical function, two-dimensional polynomial and two-dimensional trigonometrical function. It was demonstrated that this new type approximation method not only has the reproducing property of complete algebraical polynomials, but also is feasible for functions with varied boundary conditions and has excellent uniform and simultaneous convergence of functions and their (partial) derivatives up to $2r$ order. However, it is necessary to point out that the form of functions is still an important factor on the performance of this new type approximation method. For illustration, we take a comparative numerical experiment to investigate how the convergence and approximation accuracy of the algebraical polynomial interpolation based composite Fourier series method are influenced by the change of the form of functions.

In the numerical experiment, the sample functions take the form of hyperbolic sine function as shown in Eq. (1)

$$u(x_1) = \frac{\sinh[\alpha_0(a-x_1)]}{\sinh(\alpha_0 a)}, \quad x_1 \in [0,a], \tag{1}$$

where the computational parameter $\alpha_0 a$ is adjusted from 0.01 to 1.00, 10.0 and 100 successively. With the increase of computational parameter $\alpha_0 a$, the evolution of characteristics of the sample functions and approximation accuracy of the algebraical polynomial interpolation based composite Fourier series method are displayed respectively in Figures 1 and 2.

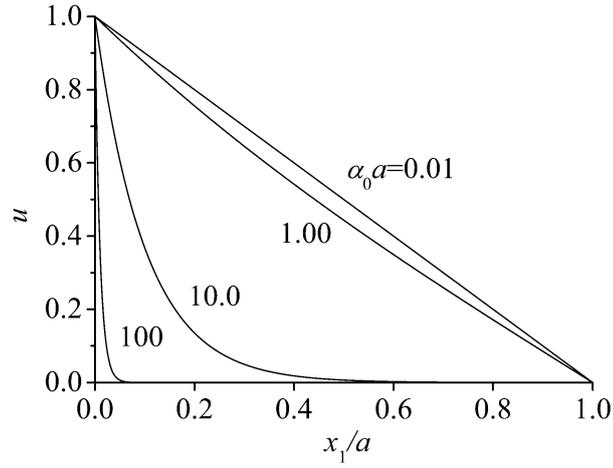

Figure 1: Sample functions with different values of computational parameter $\alpha_0 a$.

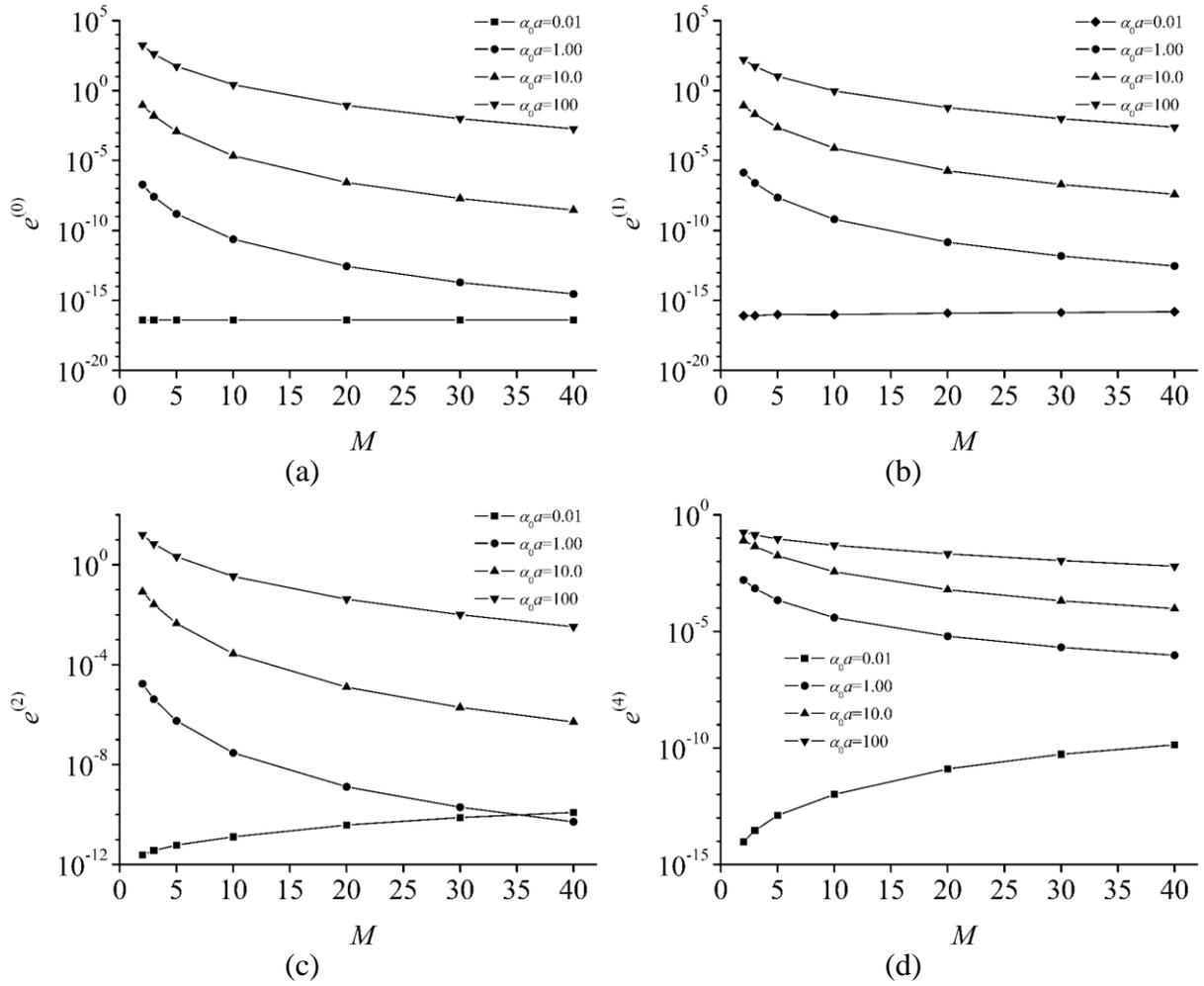

Figure 2: Convergence of the composite Fourier series for sample functions with different values of computational parameter $\alpha_0 a$: (a) $e^{(0)}(u_M)$ -$M$ curves, (b) $e^{(1)}(u_M)$ -$M$ curves, (c) $e^{(2)}(u_M)$ -$M$ curves, (d) $e^{(4)}(u_M)$ -$M$ curves.

It is observed that for the parameter $\alpha_0 a$ equaling to 0.01 or 1.00, the sample functions undergo gentle changes on the interval $[0, a]$, and the algebraical polynomial interpolation based composite Fourier series method converges rapidly and has high approximation accuracy. For the parameter $\alpha_0 a$ equaling to 10.0, the sample function slopes quickly down to zero on the interval $[0, a]$, and the algebraical polynomial interpolation based composite Fourier series method converges well and has relatively high approximation accuracy. For the parameter $\alpha_0 a$ equaling to 100, the sample function brings about the presence of the boundary layer on the interval $[0, a]$, and the algebraical polynomial interpolation based composite Fourier series method does not converge well and has low approximation accuracy.

Therefore, it is expected that the algebraical polynomial interpolation based composite Fourier series method is not appropriate for the direct application to the analytical solution of multiscale problems. Accordingly, by means of the adaptive adjustment of the boundary functions within the theoretical framework of the composite Fourier series method, we are at a starting point for the development of the new type Fourier series multiscale method.

## 3. Decomposition structures of the solution

On the basis of the decomposition structures of functions presented in [81], the Fourier series multiscale solution of the differential equations can be expressed as a superposition of the general solution and the particular solution.

### *3.1. For the full-range expansion of one-dimensional solutions*

Let $r$ be a positive integer and consider a $2r$-th order linear differential equation with constant coefficients on the interval $[-a, a]$

$$\mathcal{L} u = f, \quad x_1 \in (-a, a), \tag{2}$$

where the differential operator

$$\mathcal{L} = \sum_{k_1=0}^{2r} a_{k_1} \frac{d^{k_1}}{dx_1^{k_1}}, \tag{3}$$

and $a_{k_1}$ ($k_1 = 0, 1, \cdots, 2r$) are constant coefficients.

Taking the boundary $x_1 = a$ as an example, we have the corresponding boundary conditions

$$\mathbf{B} u = \mathbf{g}, \tag{4}$$

where the vector of the differential operators

$$\mathbf{B}^{\mathrm{T}} = [\mathcal{B}_1 \quad \mathcal{B}_2 \quad \cdots \quad \mathcal{B}_r], \tag{5}$$

and for $l = 1, 2, \cdots, r$, the differential operator

$$\mathcal{B}_l = \sum_{k_1=0}^{2r-1} b_{k_1}^l \frac{d^{k_1}}{dx_1^{k_1}}, \tag{6}$$

where $b_{k_1}^l$ are constant coefficients and $\mathbf{g}$ is a constant vector related to the prescribed boundary conditions.

Suppose that $u(x_1)$, the solution of Eq. (2), and its first to $2r$-th order derivatives can be expanded in full-range Fourier series on the interval $[-a, a]$. Then referring to the decomposition structure of functions presented in [81], we can express $u(x_1)$ as the following composite Fourier series

$$u(x_1) = \varphi_0(x_1) + \varphi_1(x_1), \tag{7}$$

where $\varphi_1(x_1)$ is the boundary function such that

$$\varphi_1^{(k_1)}(a) - \varphi_1^{(k_1)}(-a) = u^{(k_1)}(a) - u^{(k_1)}(-a), \quad k_1 = 0, 1, \cdots, 2r-1, \tag{8}$$

and $\varphi_0(x_1)$ is the internal function and naturally satisfies the sufficient conditions for $2r$ times term-by-term differentiation of the full-range Fourier series expansion of one-dimensional functions

$$\varphi_0^{(k_1)}(a) - \varphi_0^{(k_1)}(-a) = 0, \quad k_1 = 0, 1, \cdots, 2r-1. \tag{9}$$

Substitute Eq. (7) into Eq. (2), and further, suppose that

$$\mathcal{L}\varphi_1 = 0, \tag{10}$$

and

$$\mathcal{L}\varphi_0 = f, \tag{11}$$

then we can accordingly decompose the solution of Eq. (2) into two parts, namely, the general solution and the particular solution. In this decomposition, the general solution corresponds to the boundary function $\varphi_1(x_1)$ of the composite Fourier series, such that

$$\left. \begin{array}{l} \mathcal{L}\varphi_1 = 0 \\ \varphi_1^{(k_1)}(a) - \varphi_1^{(k_1)}(-a) = u^{(k_1)}(a) - u^{(k_1)}(-a), \quad k_1 = 0, 1, \cdots, 2r-1 \end{array} \right\}, \tag{12}$$

and the particular solution corresponds to the internal function $\varphi_0(x_1)$ of the composite Fourier series, such that

$$\left. \begin{array}{l} \mathcal{L}\varphi_0 = f \\ \varphi_0^{(k_1)}(a) - \varphi_0^{(k_1)}(-a) = 0, \quad k_1 = 0, 1, \cdots, 2r-1 \end{array} \right\}. \tag{13}$$

It is necessary to point out that the general solution and the particular solution of Eq. (2) are combined to satisfy the prescribed boundary conditions in Eq. (4).

*3.2. For the half-range cosine expansion of one-dimensional solutions*

Let $r$ be a positive integer and consider a $2r$-th order linear differential equation with constant coefficients on the interval $[0, a]$

$$\mathcal{L}u = f, \quad x_1 \in (0, a), \tag{14}$$

where the differential operator

$$\mathcal{L} = \sum_{k_1=0}^{r} a_{2k_1} \frac{d^{2k_1}}{dx_1^{2k_1}}, \tag{15}$$

and $a_{2k_1}$ ($k_1 = 0, 1, \cdots, r$) are constant coefficients.

Taking the boundary $x_1 = a$ as an example, we have the typical boundary conditions

$$\mathbf{B}u = \mathbf{g}, \tag{16}$$

where the vector of the differential operators

$$\mathbf{B}^{\mathrm{T}} = [\mathcal{B}_1 \quad \mathcal{B}_2 \quad \cdots \quad \mathcal{B}_r], \tag{17}$$

and for $l = 1, 2, \cdots, r$, the differential operator

$$\mathcal{B}_l = \sum_{k_1=0}^{r-1} b_{2k_1}^l \frac{d^{2k_1}}{dx_1^{2k_1}}, \tag{18}$$

or

$$\mathcal{B}_l = \sum_{k_1=0}^{r-1} b_{2k_1+1}^l \frac{d^{2k_1+1}}{dx_1^{2k_1+1}}, \tag{19}$$

where $b_{2k_1}^l$ and $b_{2k_1+1}^l$ are constant coefficients and **g** is a constant vector related to the prescribed boundary conditions.

Suppose that $u(x_1)$, the solution of Eq. (14), can be expanded in half-range cosine series on the interval $[0,a]$, and its first to $2r$-th order derivatives can be expanded in half-range Fourier series on the interval $[0,a]$. Then, by the similar procedure used in section 3.1, we can naturally decompose $u(x_1)$ into the two parts, the general solution and the particular solution. In this decomposition, the general solution corresponds to the boundary function $\varphi_1(x_1)$ of the composite Fourier series, such that

$$\left.\begin{array}{l}\mathcal{L}\varphi_1 = 0 \\ \varphi_1^{(2k_1+1)}(a) = u^{(2k_1+1)}(a), \quad \varphi_1^{(2k_1+1)}(0) = u^{(2k_1+1)}(0), \quad k_1 = 0, 1, \cdots, r-1\end{array}\right\}, \tag{20}$$

and the particular solution corresponds to the internal function $\varphi_0(x_1)$ of the composite Fourier series, such that

$$\left.\begin{array}{l}\mathcal{L}\varphi_0 = f \\ \varphi_0^{(2k_1+1)}(a) = 0, \quad \varphi_0^{(2k_1+1)}(0) = 0, \quad k_1 = 0, 1, \cdots, r-1\end{array}\right\}. \tag{21}$$

It is necessary to point out that the general solution and the particular solution of Eq. (14) are combined to satisfy the prescribed boundary conditions in Eq. (16).

*3.3. For the half-range sine expansion of one-dimensional solutions*

Let $r$ be a positive integer and let $u(x_1)$ be the solution of Eq. (14), which is a $2r$-th order linear differential equation with constant coefficients on the interval $[0,a]$. Suppose that $u(x_1)$ can be expanded in half-range sine series on the interval $[0,a]$, and its first to $2r$-th order derivatives can be expanded in half-range Fourier series on the interval $[0,a]$. Then, by the similar procedure used in section 3.1, we can naturally decompose $u(x_1)$ into the two parts, the general solution and the particular solution. In this decomposition, the general solution corresponds to the boundary function $\varphi_1(x_1)$ of the composite Fourier series, such that

$$\left.\begin{array}{l}\mathcal{L}\varphi_1 = 0 \\ \varphi_1^{(2k_1)}(a) = u^{(2k_1)}(a), \quad \varphi_1^{(2k_1)}(0) = u^{(2k_1)}(0), \quad k_1 = 0, 1, \cdots, r-1\end{array}\right\}, \tag{22}$$

and the particular solution corresponds to the internal function $\varphi_0(x_1)$ of the composite Fourier series, such that

$$\left.\begin{array}{l}\mathcal{L}\varphi_0 = f \\ \varphi_0^{(2k_1)}(a) = 0, \quad \varphi_0^{(2k_1)}(0) = 0, \quad k_1 = 0, 1, \cdots, r-1\end{array}\right\}. \tag{23}$$

It is necessary to point out that the general solution and the particular solution of Eq. (14) are combined to satisfy the prescribed boundary conditions in Eq. (16).

*3.4. For the full-range expansion of two-dimensional solutions*

Let $r$ be a positive integer and consider a $2r$-th order linear differential equation with constant coefficients in the domain $[-a,a]\times[-b,b]$

$$\mathcal{L}u = f, \quad (x_1, x_2) \in (-a,a)\times(-b,b), \tag{24}$$

where the differential operator

$$\mathcal{L} = \sum_{\substack{k_1 \geq 0, k_2 \geq 0 \\ k_1 + k_2 \leq 2r}} a_{k_1, k_2} \frac{\partial^{k_1+k_2}}{\partial x_1^{k_1} \partial x_2^{k_2}}, \tag{25}$$

and $a_{k_1, k_2}$ ($k_1, k_2 = 0, 1, 2, \cdots, k_1 + k_2 \leq 2r$) are constant coefficients.

Taking the boundary $x_1 = a$ as an example, we have the corresponding boundary conditions

$$\mathbf{B}u = \mathbf{g}, \tag{26}$$

where the vector of the differential operators

$$\mathbf{B}^{\mathrm{T}} = [\mathcal{B}_1 \quad \mathcal{B}_2 \quad \cdots \quad \mathcal{B}_r], \tag{27}$$

and for $l = 1, 2, \cdots, r$, the differential operator

$$\mathcal{B}_l = \sum_{\substack{k_1 \geq 0, k_2 \geq 0 \\ k_1 + k_2 \leq 2r-1}} b^l_{k_1, k_2} \frac{\partial^{k_1+k_2}}{\partial x_1^{k_1} \partial x_2^{k_2}}, \tag{28}$$

where $b^l_{k_1, k_2}$ are constant coefficients and $\mathbf{g}$ is a vector of functions related to the prescribed boundary conditions.

Suppose that $u(x_1, x_2)$, the solution of Eq. (24), and its first to $2r$-th order partial derivatives can be expanded in full-range Fourier series in the domain $[-a,a]\times[-b,b]$. Then referring to the decomposition structure of functions presented in [81], we can express $u(x_1, x_2)$ as the following composite Fourier series

$$u(x_1, x_2) = \varphi_0(x_1, x_2) + \varphi_1(x_1, x_2) + \varphi_2(x_1, x_2) + \varphi_3(x_1, x_2), \tag{29}$$

where $\varphi_3(x_1, x_2)$ is the corner function such that

$$\varphi_3^{(k_1,k_2)}(a,b) - \varphi_3^{(k_1,k_2)}(a,-b) - \varphi_3^{(k_1,k_2)}(-a,b) + \varphi_3^{(k_1,k_2)}(-a,-b) =$$

$$u^{(k_1,k_2)}(a,b) - u^{(k_1,k_2)}(a,-b) - u^{(k_1,k_2)}(-a,b) + u^{(k_1,k_2)}(-a,-b),$$

$$k_1, k_2 = 0, 1, 2, \cdots, \quad k_1 + k_2 \leq 2r-2, \tag{30}$$

$\varphi_1(x_1, x_2)$ and $\varphi_2(x_1, x_2)$ are the boundary functions, and respectively, satisfy the following conditions

$$\varphi_1^{(k_1,0)}(a,x_2)-\varphi_1^{(k_1,0)}(-a,x_2)=[u^{(k_1,0)}(a,x_2)-\varphi_3^{(k_1,0)}(a,x_2)]$$

$$-[u^{(k_1,0)}(-a,x_2)-\varphi_3^{(k_1,0)}(-a,x_2)],\ x_2\in(-b,b),\ k_1=0,1,\cdots,2r-1, \quad (31)$$

$$\varphi_1^{(0,k_2)}(x_1,b)-\varphi_1^{(0,k_2)}(x_1,-b)=0,\ x_1\in(-a,a),\ k_2=0,1,\cdots,2r-1, \quad (32)$$

and

$$\varphi_2^{(k_1,0)}(a,x_2)-\varphi_2^{(k_1,0)}(-a,x_2)=0,\ x_2\in(-b,b),\ k_1=0,1,\cdots,2r-1, \quad (33)$$

$$\varphi_2^{(0,k_2)}(x_1,b)-\varphi_2^{(0,k_2)}(x_1,-b)=[u^{(0,k_2)}(x_1,b)-\varphi_3^{(0,k_2)}(x_1,b)]$$

$$-[u^{(0,k_2)}(x_1,-b)-\varphi_3^{(0,k_2)}(x_1,-b)],\ x_1\in(-a,a),\ k_2=0,1,\cdots,2r-1, \quad (34)$$

and $\varphi_0(x_1,x_2)$ is the internal function and naturally satisfies the sufficient conditions for $2r$ times term-by-term differentiation of the full-range Fourier series expansion of two-dimensional functions

$$\varphi_0^{(k_1,0)}(a,x_2)-\varphi_0^{(k_1,0)}(-a,x_2)=0,\ x_2\in(-b,b),\ k_1=0,1,\cdots,2r-1, \quad (35)$$

$$\varphi_0^{(0,k_2)}(x_1,b)-\varphi_0^{(0,k_2)}(x_1,-b)=0,\ x_1\in(-a,a),\ k_2=0,1,\cdots,2r-1, \quad (36)$$

$$\varphi_0^{(k_1,k_2)}(a,b)-\varphi_0^{(k_1,k_2)}(a,-b)-\varphi_0^{(k_1,k_2)}(-a,b)+\varphi_0^{(k_1,k_2)}(-a,-b)=0,$$

$$k_1,k_2=0,1,2,\cdots,\ k_1+k_2\le 2r-2. \quad (37)$$

Substitute Eq. (29) into Eq. (24), and further, suppose that

$$\mathcal{L}\varphi_1=0, \quad (38)$$

$$\mathcal{L}\varphi_2=0, \quad (39)$$

and

$$\mathcal{L}(\varphi_0+\varphi_3)=f, \quad (40)$$

then we can correspondingly decompose the solution of Eq. (24) into the two parts, the general solution and the particular solution. In this decomposition, the general solution corresponds to the boundary functions $\varphi_1(x_1,x_2)$ and $\varphi_2(x_1,x_2)$ that respectively satisfy the following equations

$$\left.\begin{array}{l}\mathcal{L}\varphi_1=0\\ \varphi_1^{(k_1,0)}(a,x_2)-\varphi_1^{(k_1,0)}(-a,x_2)=[u^{(k_1,0)}(a,x_2)-\varphi_3^{(k_1,0)}(a,x_2)]\\ -[u^{(k_1,0)}(-a,x_2)-\varphi_3^{(k_1,0)}(-a,x_2)],\ x_2\in(-b,b),\ k_1=0,1,\cdots,2r-1\\ \varphi_1^{(0,k_2)}(x_1,b)-\varphi_1^{(0,k_2)}(x_1,-b)=0,\ x_1\in(-a,a),\ k_2=0,1,\cdots,2r-1\end{array}\right\}, \quad (41)$$

and

$$\left.\begin{aligned}
&\mathcal{L}\varphi_2 = 0 \\
&\varphi_2^{(k_1,0)}(a,x_2) - \varphi_2^{(k_1,0)}(-a,x_2) = 0, \quad x_2 \in (-b,b), \quad k_1 = 0, 1, \cdots, 2r-1 \\
&\varphi_2^{(0,k_2)}(x_1,b) - \varphi_2^{(0,k_2)}(x_1,-b) = [u^{(0,k_2)}(x_1,b) - \varphi_3^{(0,k_2)}(x_1,b)] \\
&\quad - [u^{(0,k_2)}(x_1,-b) - \varphi_3^{(0,k_2)}(x_1,-b)], \quad x_1 \in (-a,a), \quad k_2 = 0, 1, \cdots, 2r-1
\end{aligned}\right\}. \quad (42)$$

The particular solution corresponds to the sum of the internal function $\varphi_0(x_1,x_2)$ and the corner function $\varphi_3(x_1,x_2)$ of the composite Fourier series, such that

$$\left.\begin{aligned}
&\mathcal{L}(\varphi_0 + \varphi_3) = f \\
&\varphi_3^{(k_1,k_2)}(a,b) - \varphi_3^{(k_1,k_2)}(a,-b) - \varphi_3^{(k_1,k_2)}(-a,b) + \varphi_3^{(k_1,k_2)}(-a,-b) = \\
&u^{(k_1,k_2)}(a,b) - u^{(k_1,k_2)}(a,-b) - u^{(k_1,k_2)}(-a,b) + u^{(k_1,k_2)}(-a,-b), \\
&k_1, k_2 = 0, 1, 2, \cdots, \quad k_1 + k_2 \leq 2r-2
\end{aligned}\right\}. \quad (43)$$

It is necessary to point out that the general solution and the particular solution of Eq. (24) are combined to satisfy the prescribed boundary conditions in Eq. (26).

*3.5. For the half-range sine-sine expansion of two-dimensional solutions*

Let $r$ be a positive integer and consider a $2r$-th order linear differential equation with constant coefficients in the domain $[0,a] \times [0,b]$

$$\mathcal{L}u = f, \quad (x_1, x_2) \in (0,a) \times (0,b), \quad (44)$$

where the differential operator

$$\mathcal{L} = \sum_{\substack{k_1 \geq 0, k_2 \geq 0 \\ k_1 + k_2 \leq r}} a_{2k_1, 2k_2} \frac{\partial^{2k_1 + 2k_2}}{\partial x_1^{2k_1} \partial x_2^{2k_2}}, \quad (45)$$

and $a_{2k_1, 2k_2}$ $(k_1, k_2 = 0, 1, 2, \cdots, k_1 + k_2 \leq r)$ are constant coefficients.

Taking the boundary $x_1 = a$ as an example, we have the typical boundary conditions

$$\mathbf{B}u = \mathbf{g}, \quad (46)$$

where the vector of the differential operators

$$\mathbf{B}^T = [\mathcal{B}_1 \quad \mathcal{B}_2 \quad \cdots \quad \mathcal{B}_r], \quad (47)$$

and for $l = 1, 2, \cdots, r$, the differential operator

$$\mathcal{B}_l = \sum_{\substack{k_1 \geq 0, k_2 \geq 0 \\ k_1 + k_2 \leq r-1}} b_{2k_1, 2k_2}^l \frac{\partial^{2k_1 + 2k_2}}{\partial x^{2k_1} \partial y^{2k_2}}, \quad (48)$$

or

$$\mathcal{B}_l = \sum_{\substack{k_1 \geq 0, k_2 \geq 0 \\ k_1 + k_2 \leq r-1}} b_{2k_1+1, 2k_2}^l \frac{\partial^{2k_1 + 2k_2 + 1}}{\partial x^{2k_1+1} \partial y^{2k_2}}, \quad (49)$$

or

$$\mathcal{B}_l = \sum_{\substack{k_1 \geq 0, k_2 \geq 0 \\ k_1 + k_2 \leq r-1}} b_{2k_1, 2k_2+1}^l \frac{\partial^{2k_1 + 2k_2 + 1}}{\partial x^{2k_1} \partial y^{2k_2+1}}, \quad (50)$$

or

$$\mathcal{B}_l = \sum_{\substack{k_1 \geq 0, k_2 \geq 0 \\ k_1+k_2 \leq r-2}} b^l_{2k_1+1, 2k_2+1} \frac{\partial^{2k_1+2k_2+2}}{\partial x^{2k_1+1} \partial y^{2k_2+1}}, \tag{51}$$

where $b^l_{2k_1,2k_2}$, $b^l_{2k_1+1,2k_2}$, $b^l_{2k_1,2k_2+1}$ and $b^l_{2k_1+1,2k_2+1}$ are constant coefficients and **g** is a vector of functions related to the prescribed boundary conditions.

Suppose that $u(x_1, x_2)$, the solution of Eq. (44), can be expanded in half-range sine-sine series in the domain $[0,a] \times [0,b]$, and its first to $2r$-th order partial derivatives can be expanded in half-range Fourier series in the domain $[0,a] \times [0,b]$. Then, by the similar procedure used in section 3.4, we can naturally decompose $u(x_1, x_2)$ into the two parts, the general solution and the particular solution. In this decomposition, the general solution corresponds to the boundary functions $\varphi_1(x_1, x_2)$ and $\varphi_2(x_1, x_2)$ that respectively satisfy the following equations

$$\left.\begin{aligned}
& \mathcal{L}\varphi_1 = 0 \\
& \varphi_1^{(2k_1,0)}(a, x_2) = u^{(2k_1,0)}(a, x_2) - \varphi_3^{(2k_1,0)}(a, x_2), \\
& \varphi_1^{(2k_1,0)}(0, x_2) = u^{(2k_1,0)}(0, x_2) - \varphi_3^{(2k_1,0)}(0, x_2), \quad x_2 \in (0,b), \quad k_1 = 0, 1, \cdots, r-1 \\
& \varphi_1^{(0,2k_2)}(x_1, b) = 0, \quad \varphi_1^{(0,2k_2)}(x_1, 0) = 0, \quad x_1 \in (0,a), \quad k_2 = 0, 1, \cdots, r-1
\end{aligned}\right\}, \tag{52}$$

and

$$\left.\begin{aligned}
& \mathcal{L}\varphi_2 = 0 \\
& \varphi_2^{(2k_1,0)}(a, x_2) = 0, \quad \varphi_2^{(2k_1,0)}(0, x_2) = 0, \quad x_2 \in (0,b), \quad k_1 = 0, 1, \cdots, r-1 \\
& \varphi_2^{(0,2k_2)}(x_1, b) = u^{(0,2k_2)}(x_1, b) - \varphi_3^{(0,2k_2)}(x_1, b), \\
& \varphi_2^{(0,2k_2)}(x_1, 0) = u^{(0,2k_2)}(x_1, 0) - \varphi_3^{(0,2k_2)}(x_1, 0), \quad x_1 \in (0,a), \quad k_2 = 0, 1, \cdots, r-1
\end{aligned}\right\}. \tag{53}$$

The particular solution corresponds to the sum of the internal function $\varphi_0(x_1, x_2)$ and the corner function $\varphi_3(x_1, x_2)$ of the composite Fourier series, such that

$$\left.\begin{aligned}
& \mathcal{L}(\varphi_0 + \varphi_3) = f \\
& \varphi_3^{(2k_1,2k_2)}(a,b) = u^{(2k_1,2k_2)}(a,b), \quad \varphi_3^{(2k_1,2k_2)}(a,0) = u^{(2k_1,2k_2)}(a,0), \\
& \varphi_3^{(2k_1,2k_2)}(0,b) = u^{(2k_1,2k_2)}(0,b), \quad \varphi_3^{(2k_1,2k_2)}(0,0) = u^{(2k_1,2k_2)}(0,0), \\
& k_1, k_2 = 0, 1, 2, \cdots, \quad k_1 + k_2 \leq r-1
\end{aligned}\right\}. \tag{54}$$

It is necessary to point out that the general solution and the particular solution of Eq. (44) are combined to satisfy the prescribed boundary conditions in Eq. (46).

## 4. Derivation of the general solution

As shown in section 3, the general solution of the linear differential equation with constant coefficients is required not only to satisfy the selection criteria of boundary functions, but also to satisfy the corresponding homogeneous form of the original differential equation. In this section, we present the derivation method of the general solution with examples of one-dimensional full-range expansion and two-dimensional full-range expansion.

*4.1. For the full-range expansion of one-dimensional solutions*

For the full-range expansion of the one-dimensional solution, without loss of generality,

we suppose that Eq. (12) has the following homogeneous solutions in exponential form
$$p_H(x_1) = \exp(\eta x_1), \tag{55}$$
where $\eta$ is an undetermined constant.

Substitute Eq. (55) into Eq. (12), then we can obtain the characteristic equation
$$\sum_{k_1=0}^{2r} a_{k_1} \eta^{k_1} = 0, \tag{56}$$
and the corresponding $2r$ characteristic roots $\{\eta_l\}_{1 \leq l \leq 2r}$.

Further, according to the distribution of all characteristic roots of the characteristic equation, we can determine the $2r$ linearly independent homogeneous solutions $\{p_{l,H}(x_1)\}_{1 \leq l \leq 2r}$ of Eq. (12).

We select the vector of functions
$$\mathbf{p}_1^T(x_1) = [p_{1,H}(x_1) \quad p_{2,H}(x_1) \quad \cdots \quad p_{2r,H}(x_1)], \tag{57}$$
and substitute it into Eq. (15) of [81], then we can express the general solution as
$$\varphi_1(x_1) = \mathbf{\Phi}_1^T(x_1) \cdot \mathbf{q}_1, \tag{58}$$
where the definitions of the vector of basis function $\mathbf{\Phi}_1^T(x_1)$ and the vector of undetermined constants $\mathbf{q}_1^T$ are respectively referred to Eqs. (22) and (20) in [81].

*4.2. For the full-range expansion of two-dimensional solutions*

For the full-range expansion of the two-dimensional solution, we only illustrate the procedures for constructing the general solution $\varphi_1(x_1, x_2)$.

For this purpose, we expand the function $\varphi_1(x_1, x_2)$ in a one-dimensional full-range Fourier series on the interval $[-b, b]$ along the $x_2$ direction, and then we have
$$\varphi_1(x_1, x_2) = \sum_{n=0}^{\infty} \mu_n [\xi_{1n}(x_1) \cos(\beta_n x_2) + \xi_{2n}(x_1) \sin(\beta_n x_2)], \tag{59}$$
where $\beta_n = n\pi/b$, $\mu_n = \begin{cases} 1/2, & n = 0 \\ 1, & n > 0 \end{cases}$, $\xi_{1n}(x_1)$ and $\xi_{2n}(x_1)$ are the corresponding one-dimensional full-range Fourier coefficients.

Substitute it into Eq. (41), then we find

1. for $n > 0$, the equations of the undetermined functions $\xi_{1n}(x_1)$ and $\xi_{2n}(x_1)$ are
$$\begin{bmatrix} \mathcal{L}_{1n,1} & \mathcal{L}_{1n,2} \\ \mathcal{L}_{2n,1} & \mathcal{L}_{2n,2} \end{bmatrix} \begin{bmatrix} \xi_{1n} \\ \xi_{2n} \end{bmatrix} = \mathbf{0}, \tag{60}$$

where the differential operators
$$\mathcal{L}_{1n,1} = \mathcal{L}_{2n,2} = \sum_{\substack{k_1 \geq 0, k_2 \geq 0 \\ k_1 + k_2 \leq r}} a_{2k_1, 2k_2} (-1)^{k_2} \beta_n^{2k_2} \frac{d^{2k_1}}{dx_1^{2k_1}} + \sum_{\substack{k_1 \geq 0, k_2 \geq 0 \\ k_1 + k_2 \leq r-1}} a_{2k_1+1, 2k_2} (-1)^{k_2} \beta_n^{2k_2} \frac{d^{2k_1+1}}{dx_1^{2k_1+1}}, \tag{61}$$

$$\mathcal{L}_{1n,2} = \sum_{\substack{k_1 \geq 0, k_2 \geq 0 \\ k_1 + k_2 \leq r-1}} a_{2k_1, 2k_2+1} (-1)^{k_2} \beta_n^{2k_2+1} \frac{d^{2k_1}}{dx_1^{2k_1}} + \sum_{\substack{k_1 \geq 0, k_2 \geq 0 \\ k_1 + k_2 \leq r-1}} a_{2k_1+1, 2k_2+1} (-1)^{k_2} \beta_n^{2k_2+1} \frac{d^{2k_1+1}}{dx_1^{2k_1+1}}, \tag{62}$$

$$\mathcal{L}_{2n,1} = \sum_{\substack{k_1 \geq 0, k_2 \geq 0 \\ k_1+k_2 \leq r-1}} a_{2k_1,2k_2+1}(-1)^{k_2+1} \beta_n^{2k_2+1} \frac{d^{2k_1}}{dx_1^{2k_1}} + \sum_{\substack{k_1 \geq 0, k_2 \geq 0 \\ k_1+k_2 \leq r-1}} a_{2k_1+1,2k_2+1}(-1)^{k_2+1} \beta_n^{2k_2+1} \frac{d^{2k_1+1}}{dx_1^{2k_1+1}}.$$

(63)

2. for $n = 0$, the equation of the undetermined function $\xi_{10}(x_1)$ is

$$\mathcal{L}_{10}\xi_{10} = 0,$$ (64)

where the differential operator

$$\mathcal{L}_{10} = \sum_{0 \leq k_1 \leq 2r} a_{k_1,0} \frac{d^{k_1}}{dx_1^{k_1}}.$$ (65)

1. Solving Eq. (60)

Further, let the solutions of Eq. (60) be in the exponential form

$$\left. \begin{array}{l} \xi_{1n}(x_1) = G_{1n,1} p_{1n}(x_1) \\ \xi_{2n}(x_1) = G_{1n,2} p_{1n}(x_1) \end{array} \right\},$$ (66)

where the basis function

$$p_{1n}(x_1) = \exp(\eta_n x_1),$$ (67)

$G_{1n,1}$, $G_{1n,2}$ and $\eta_n$ are undetermined constants.

Substitute the equation above into Eq. (60), then we have

$$\begin{bmatrix} t_{1n,1} & t_{1n,2} \\ t_{2n,1} & t_{2n,2} \end{bmatrix} \begin{bmatrix} G_{1n,1} \\ G_{1n,2} \end{bmatrix} = 0,$$ (68)

where

$$t_{1n,1} = t_{2n,2} = \sum_{\substack{k_1 \geq 0, k_2 \geq 0 \\ k_1+k_2 \leq r}} a_{2k_1,2k_2}(-1)^{k_2} \eta_n^{2k_1} \beta_n^{2k_2} + \sum_{\substack{k_1 \geq 0, k_2 \geq 0 \\ k_1+k_2 \leq r-1}} a_{2k_1+1,2k_2}(-1)^{k_2} \eta_n^{2k_1+1} \beta_n^{2k_2},$$ (69)

$$t_{1n,2} = \sum_{\substack{k_1 \geq 0, k_2 \geq 0 \\ k_1+k_2 \leq r-1}} a_{2k_1,2k_2+1}(-1)^{k_2} \eta_n^{2k_1} \beta_n^{2k_2+1} + \sum_{\substack{k_1 \geq 0, k_2 \geq 0 \\ k_1+k_2 \leq r-1}} a_{2k_1+1,2k_2+1}(-1)^{k_2} \eta_n^{2k_1+1} \beta_n^{2k_2+1},$$ (70)

$$t_{2n,1} = \sum_{\substack{k_1 \geq 0, k_2 \geq 0 \\ k_1+k_2 \leq r-1}} a_{2k_1,2k_2+1}(-1)^{k_2+1} \eta_n^{2k_1} \beta_n^{2k_2+1} + \sum_{\substack{k_1 \geq 0, k_2 \geq 0 \\ k_1+k_2 \leq r-1}} a_{2k_1+1,2k_2+1}(-1)^{k_2+1} \eta_n^{2k_1+1} \beta_n^{2k_2+1}.$$ (71)

As to the characteristic equation of Eq. (68), we can obtain $4r$ characteristic roots $\{\eta_{n,l}\}_{1 \leq l \leq 4r}$ and the linearly independent homogeneous solutions $\{p_{1n,l}(x_1)\}_{1 \leq l \leq 4r}$.

For $l = 1, 2, \cdots, 4r$, we substitute the characteristic roots $\eta_{n,l}$ into Eq. (68), and further, obtain the following relations between the undetermined coefficients $G_{1n,1}^l$ and $G_{1n,2}^l$

$$\begin{bmatrix} t_{1n,1}(\eta_{n,l}) & t_{1n,2}(\eta_{n,l}) \end{bmatrix} \begin{bmatrix} G_{1n,1}^l \\ G_{1n,2}^l \end{bmatrix} = 0.$$ (72)

If we write the undetermined constant vectors

$$\mathbf{a}_{1n,1}^{\mathrm{T}} = [G_{1n,1}^1 \quad G_{1n,1}^2 \quad \cdots \quad G_{1n,1}^{4r}],$$ (73)

$$\mathbf{a}_{1n,2}^{\mathrm{T}} = [G_{1n,2}^1 \quad G_{1n,2}^2 \quad \cdots \quad G_{1n,2}^{4r}],$$ (74)

then the above equations can be combined and expressed as

$$\begin{bmatrix} \mathbf{a}_{1n,1} \\ \mathbf{a}_{1n,2} \end{bmatrix} = \mathbf{T}_{1n} \mathbf{a}_{1n,1},$$ (75)

where $\mathbf{T}_{1n}$ is a transformation matrix.

Referring to section 3.2 in [81], we can select a vector of functions
$$\mathbf{p}_{1n}^{\mathrm{T}}(x_1) = [p_{1n,1}(x_1) \quad p_{1n,2}(x_1) \quad \cdots \quad p_{1n,4r}(x_1)], \tag{76}$$
and construct the two functions $\xi_{1n}(x_1)$ and $\xi_{2n}(x_1)$
$$\xi_{1n}(x_1) = \mathbf{p}_{1n}^{\mathrm{T}}(x_1) \cdot \mathbf{a}_{1n,1}, \tag{77}$$
$$\xi_{2n}(x_1) = \mathbf{p}_{1n}^{\mathrm{T}}(x_1) \cdot \mathbf{a}_{1n,2}. \tag{78}$$
Meanwhile, referring to Eqs. (54) and (55) in [81], we have
$$\mathbf{R}_{1n}\mathbf{a}_{1n,1} = \mathbf{q}_{1n,1}, \tag{79}$$
$$\mathbf{R}_{1n}\mathbf{a}_{1n,2} = \mathbf{q}_{1n,2}, \tag{80}$$
where the definitions of the matrix $\mathbf{R}_{1n}$, the sub-vectors of boundary Fourier coefficients $\mathbf{q}_{1n,1}$ and $\mathbf{q}_{1n,2}$ are referred to Eqs. (64)-(66) in [81].

Further, we can rewrite Eqs. (79) and (80) and obtain
$$\mathbf{S}_{1n}\begin{bmatrix}\mathbf{a}_{1n,1}\\ \mathbf{a}_{1n,2}\end{bmatrix} = \begin{bmatrix}\mathbf{q}_{1n,1}\\ \mathbf{q}_{1n,2}\end{bmatrix}, \tag{81}$$
where the matrix
$$\mathbf{S}_{1n} = \begin{bmatrix}\mathbf{R}_{1n} & \mathbf{0}\\ \mathbf{0} & \mathbf{R}_{1n}\end{bmatrix}. \tag{82}$$
Combining Eqs. (75) and (81), we have
$$\mathbf{S}_{1n}\mathbf{T}_{1n}\mathbf{a}_{1n,1} = \begin{bmatrix}\mathbf{q}_{1n,1}\\ \mathbf{q}_{1n,2}\end{bmatrix}, \tag{83}$$
hence
$$\mathbf{a}_{1n,1} = (\mathbf{S}_{1n}\mathbf{T}_{1n})^{-1}\begin{bmatrix}\mathbf{q}_{1n,1}\\ \mathbf{q}_{1n,2}\end{bmatrix}. \tag{84}$$
Therefore, the undetermined one-dimensional functions $\xi_{1n}(x_1)$ and $\xi_{2n}(x_1)$ can be expressed as
$$\begin{bmatrix}\xi_{1n}(x_1)\\ \xi_{2n}(x_1)\end{bmatrix} = \mathbf{p}_{1n,R}^{\mathrm{T}}(x_1) \cdot \begin{bmatrix}\mathbf{q}_{1n,1}\\ \mathbf{q}_{1n,2}\end{bmatrix}, \tag{85}$$
where the vector of functions
$$\mathbf{p}_{1n,R}^{\mathrm{T}}(x_1) = \begin{bmatrix}\mathbf{p}_{1n}^{\mathrm{T}}(x_1) & \mathbf{0}\\ \mathbf{0} & \mathbf{p}_{1n}^{\mathrm{T}}(x_1)\end{bmatrix} \cdot \mathbf{T}_{1n}(\mathbf{S}_{1n}\mathbf{T}_{1n})^{-1}. \tag{86}$$
Accordingly, we can define the vector of basis functions
$$\mathbf{\Phi}_{1n}^{\mathrm{T}}(x_1, x_2) = \begin{bmatrix}\mathbf{\Phi}_{1n,1}^{\mathrm{T}}(x_1, x_2) & \mathbf{\Phi}_{1n,2}^{\mathrm{T}}(x_1, x_2)\end{bmatrix} = \mathbf{H}_{1n}(x_2) \cdot \mathbf{p}_{1n,R}^{\mathrm{T}}(x_1), \tag{87}$$
where the vector of functions
$$\mathbf{H}_{1n}(x_2) = [\cos(\beta_n x_2) \quad \sin(\beta_n x_2)], \tag{88}$$
and the corresponding vector of higher order partial derivatives of the basis functions
$$\mathbf{\Phi}_{1n}^{(k_1,k_2)\mathrm{T}}(x_1, x_2) = \begin{bmatrix}\mathbf{\Phi}_{1n,1}^{(k_1,k_2)\mathrm{T}}(x_1, x_2) & \mathbf{\Phi}_{1n,2}^{(k_1,k_2)\mathrm{T}}(x_1, x_2)\end{bmatrix} = \mathbf{H}_{1n}^{(k_2)}(x_2) \cdot \mathbf{p}_{1n,R}^{(k_1)\mathrm{T}}(x_1), \tag{89}$$
where $k_1$ and $k_2$ are nonnegative integers, and the function matrix
$$\mathbf{p}_{1n,R}^{(k_1)\mathrm{T}}(x_1) = \begin{bmatrix}\mathbf{p}_{1n}^{(k_1)\mathrm{T}}(x_1) & \mathbf{0}\\ \mathbf{0} & \mathbf{p}_{1n}^{(k_1)\mathrm{T}}(x_1)\end{bmatrix} \cdot \mathbf{T}_{1n}(\mathbf{S}_{1n}\mathbf{T}_{1n})^{-1}, \tag{90}$$
and the vector of functions

$$\mathbf{H}_{1n}^{(k_2)}(x_2) = \left[ [\cos(\beta_n x_2)]^{(k_2)} \quad [\sin(\beta_n x_2)]^{(k_2)} \right]. \tag{91}$$

It should be pointed out that the foregoing derivation method is based on the proper selection of the vector of functions, $\mathbf{p}_{1n}^{T}(x_1)$. In specific situations, if we can directly obtain the linearly independent homogeneous solutions $\{p_{1nl,H}(x_1,x_2)\}_{1 \leq l \leq 4r}$ of Eq. (41), then we construct function

$$\xi_{1n}(x_1)\cos(\beta_n x_2) + \xi_{2n}(x_1)\sin(\beta_n x_2) = \mathbf{p}_{1n,H}^{T}(x_1,x_2) \cdot \mathbf{a}_{1n}, \tag{92}$$

where the vector of functions

$$\mathbf{p}_{1n,H}^{T}(x_1,x_2) = [p_{1n1,H}(x_1,x_2) \quad p_{1n2,H}(x_1,x_2) \quad \cdots \quad p_{1n(4r),H}(x_1,x_2)], \tag{93}$$

and the undetermined constant vector

$$\mathbf{a}_{1n}^{T} = [G_{1n}^{1} \quad G_{1n}^{2} \quad \cdots \quad G_{1n}^{4r}]. \tag{94}$$

Therefore, referring to Eqs. (54)-(55) in [81], we obtain

$$\mathbf{R}_{1n,H}\mathbf{a}_{1n} = \begin{bmatrix} \mathbf{q}_{1n,1} \\ \mathbf{q}_{1n,2} \end{bmatrix}, \tag{95}$$

where the matrices

$$\mathbf{R}_{1n,H} = \begin{bmatrix} \mathbf{R}_{1n1,H} \\ \mathbf{R}_{1n2,H} \end{bmatrix}, \tag{96}$$

$$\mathbf{R}_{1n1,H} = \begin{bmatrix} \mathbf{p}_{1n,H}^{(0,0)T}(a,0) - \mathbf{p}_{1n,H}^{(0,0)T}(-a,0) \\ \mathbf{p}_{1n,H}^{(1,0)T}(a,0) - \mathbf{p}_{1n,H}^{(1,0)T}(-a,0) \\ \vdots \\ \mathbf{p}_{1n,H}^{(2r-1,0)T}(a,0) - \mathbf{p}_{1n,H}^{(2r-1,0)T}(-a,0) \end{bmatrix}, \tag{97}$$

$$\mathbf{R}_{1n2,H} = \begin{bmatrix} \mathbf{p}_{1n,H}^{(0,0)T}(a,b/2n) - \mathbf{p}_{1n,H}^{(0,0)T}(-a,b/2n) \\ \mathbf{p}_{1n,H}^{(1,0)T}(a,b/2n) - \mathbf{p}_{1n,H}^{(1,0)T}(-a,b/2n) \\ \vdots \\ \mathbf{p}_{1n,H}^{(2r-1,0)T}(a,b/2n) - \mathbf{p}_{1n,H}^{(2r-1,0)T}(-a,b/2n) \end{bmatrix}, \tag{98}$$

and the definitions of the sub-vectors of boundary Fourier coefficients $\mathbf{q}_{1n,1}$ and $\mathbf{q}_{1n,2}$ are respectively referred to Eqs. (65), (66) in [81].

Accordingly, we can define the vector of basis functions

$$\mathbf{\Phi}_{1n}^{T}(x_1,x_2) = \left[\mathbf{\Phi}_{1n,1}^{T}(x_1,x_2) \quad \mathbf{\Phi}_{1n,2}^{T}(x_1,x_2)\right] = \mathbf{p}_{1n,H}^{T}(x_1,x_2) \cdot \mathbf{R}_{1n,H}^{-1}, \tag{99}$$

and the corresponding vector of higher order partial derivatives of basis functions

$$\mathbf{\Phi}_{1n}^{(k_1,k_2)T}(x_1,x_2) = \left[\mathbf{\Phi}_{1n,1}^{(k_1,k_2)T}(x_1,x_2) \quad \mathbf{\Phi}_{1n,2}^{(k_1,k_2)T}(x_1,x_2)\right] = \mathbf{p}_{1n,H}^{(k_1,k_2)T}(x_1,x_2) \cdot \mathbf{R}_{1n,H}^{-1}, \tag{100}$$

where $k_1$ and $k_2$ are nonnegative integers.

2. Solving Eq. (64)

Let the solution of Eq. (64) be in the exponential form

$$\xi_{10}(x_1) = G_{10} p_{10}(x_1), \tag{101}$$

where the basis function

$$p_{10}(x_1) = \exp(\eta_0 x_1), \tag{102}$$

$G_{10}$ and $\eta_0$ are undetermined constants.

Substitute Eq. (102) into Eq. (64), then we can obtain the characteristic equation

$$\sum_{k_1=0}^{2r} a_{k_1,0} \eta_0^{k_1} = 0, \tag{103}$$

and the corresponding characteristic roots $\{\eta_{0,l}\}_{1\leq l\leq 2r}$, together with the linearly independent homogeneous solutions $\{p_{10,l}(x_1)\}_{1\leq l\leq 2r}$.

Referring to section 3.2 in [81], we construct the function

$$\xi_{10}(x_1) = \mathbf{p}_{10}^{\mathrm{T}}(x_1)\cdot \mathbf{a}_{10}, \tag{104}$$

where the vector of functions

$$\mathbf{p}_{10}^{\mathrm{T}}(x_1) = [p_{10,1}(x_1) \quad p_{10,2}(x_1) \quad \cdots \quad p_{10,2r}(x_1)], \tag{105}$$

and the undetermined constant vector

$$\mathbf{a}_{10,1}^{\mathrm{T}} = [G_{10,1}^1 \quad G_{10,1}^2 \quad \cdots \quad G_{10,1}^{2r}]. \tag{106}$$

Therefore, referring to Eq. (62) in [81], we obtain

$$\mathbf{R}_{10}\mathbf{a}_{10,1} = \mathbf{q}_{10,1}, \tag{107}$$

where the definitions of the matrix $\mathbf{R}_{10}$, the sub-vectors of boundary Fourier coefficients $\mathbf{q}_{10,1}$ are respectively referred to Eqs. (64), (65) in [81].

Therefore, the undetermined one-dimensional function $\xi_{10}(x_1)$ can be express as

$$\xi_{10}(x_1) = \mathbf{p}_{10,R}^{\mathrm{T}}(x_1)\cdot \mathbf{q}_{10,1}, \tag{108}$$

where the vector of functions

$$\mathbf{p}_{10,R}^{\mathrm{T}}(x_1) = \mathbf{p}_{10}^{\mathrm{T}}(x_1)\cdot \mathbf{R}_{10}^{-1}. \tag{109}$$

Accordingly, we can define the sub-vector of basis functions

$$\boldsymbol{\Phi}_{10,1}^{\mathrm{T}}(x_1,x_2) = \mathbf{H}_{10}(x_2)\cdot \mathbf{p}_{10,R}^{\mathrm{T}}(x_1), \tag{110}$$

where the vector of functions

$$\mathbf{H}_{10}(x_2) = \left[\frac{1}{2}\right], \tag{111}$$

and the corresponding vector of higher order partial derivatives of basis functions

$$\boldsymbol{\Phi}_{10,1}^{(k_1,k_2)\mathrm{T}}(x_1,x_2) = \mathbf{H}_{10}^{(k_2)}(x_2)\cdot \mathbf{p}_{10,R}^{(k_1)\mathrm{T}}(x_1), \tag{112}$$

where $k_1$ and $k_2$ are nonnegative integers.

By adopting the construction methods in section 3.2, [81], for the vectors of basis functions $\boldsymbol{\Phi}_{1,1}^{\mathrm{T}}(x_1,x_2)$, $\boldsymbol{\Phi}_{1,2}^{\mathrm{T}}(x_1,x_2)$ and $\boldsymbol{\Phi}_{1}^{\mathrm{T}}(x_1,x_2)$, and the vectors of boundary Fourier coefficients $\mathbf{q}_{1,1}^{\mathrm{T}}$, $\mathbf{q}_{1,2}^{\mathrm{T}}$ and $\mathbf{q}_{1}^{\mathrm{T}}$, we can eventually express the boundary function $\varphi_1(x_1,x_2)$ as

$$\varphi_1(x_1,x_2) = \boldsymbol{\Phi}_{1}^{\mathrm{T}}(x_1,x_2)\cdot \mathbf{q}_1. \tag{113}$$

## 5. Equivalent transformation of the solution

With the technique developed in section 4, the general solution, and further the Fourier series solution appropriate for specific kind of Fourier series expansion, are to be derived. In the obtained Fourier series multiscale solution, the primary undetermined constants include the Fourier coefficients, boundary Fourier coefficients, boundary (end) values, and sometimes the corner values of the solution to be determined. This leads to difficulties in the application of variational methods, where prior satisfaction of the displacement type boundary conditions is required. Therefore, in this section we confine ourselves to the specific two-dimensional full-range expansion in the domain $[-a,a]\times[-b,b]$ and develop the equivalent transformation method for the change of primary undetermined constants in the Fourier series multiscale solution.

For this purpose, suppose that the displacement type boundary conditions on boundary $x_1 = a$ can be expressed as

$$\bar{\mathbf{u}} = \mathbf{C}u, \tag{114}$$

where the vector of the differential operators

$$\mathbf{C}^{\mathrm{T}} = [\mathcal{C}_1 \quad \mathcal{C}_2 \quad \cdots \quad \mathcal{C}_r], \tag{115}$$

and for $l = 1, 2, \cdots, r$, the differential operator

$$\mathcal{C}_l = \sum_{\substack{k_1 \geq 0, k_2 \geq 0 \\ k_1 + k_2 \leq 2r-1}} c^l_{k_1,k_2} \frac{\partial^{k_1+k_2}}{\partial x_1^{k_1} \partial x_2^{k_2}}, \tag{116}$$

where $c^l_{k_1,k_2}$ are constant coefficients and $\bar{\mathbf{u}}$ is a vector of functions corresponding to the prescribed displacement type boundary conditions.

Substitute the Fourier series multiscale solution into Eq. (114), then we obtain

$$\bar{\mathbf{u}}(a, x_2) = \mathbf{\Gamma}(a, x_2) \cdot \mathbf{q}, \tag{117}$$

where

$$\mathbf{\Gamma} = \begin{bmatrix} \mathcal{C}_1 \mathbf{\Phi}_0^{\mathrm{T}} & \mathcal{C}_1 \mathbf{\Phi}_1^{\mathrm{T}} & \mathcal{C}_1 \mathbf{\Phi}_2^{\mathrm{T}} & \mathcal{C}_1 \mathbf{\Phi}_3^{\mathrm{T}} \\ \mathcal{C}_2 \mathbf{\Phi}_0^{\mathrm{T}} & \mathcal{C}_2 \mathbf{\Phi}_1^{\mathrm{T}} & \mathcal{C}_2 \mathbf{\Phi}_2^{\mathrm{T}} & \mathcal{C}_2 \mathbf{\Phi}_3^{\mathrm{T}} \\ \vdots & \vdots & \vdots & \vdots \\ \mathcal{C}_r \mathbf{\Phi}_0^{\mathrm{T}} & \mathcal{C}_r \mathbf{\Phi}_1^{\mathrm{T}} & \mathcal{C}_r \mathbf{\Phi}_2^{\mathrm{T}} & \mathcal{C}_r \mathbf{\Phi}_3^{\mathrm{T}} \end{bmatrix}, \tag{118}$$

and the vector of primary undetermined constants

$$\mathbf{q}^{\mathrm{T}} = [\mathbf{q}_0^{\mathrm{T}} \quad \mathbf{q}_1^{\mathrm{T}} \quad \mathbf{q}_2^{\mathrm{T}} \quad \mathbf{q}_3^{\mathrm{T}}]. \tag{119}$$

Expand $\bar{\mathbf{u}}(a, x_2)$ and the elements of the first, second, … and $r$-th rows of the matrix $\mathbf{\Gamma}(a, x_2)$ in full-range Fourier series on the interval $[-b,b]$, and let $N$ be the number of truncated terms. Then by comparing the first to $N$-th Fourier coefficients on the both sides of Eq. (117), we have

$$\mathbf{q}_{b,a+} = \mathbf{R}_{b,a+} \cdot \mathbf{q}, \tag{120}$$

where the vector of undetermined constants

$$\mathbf{q}_{b,a+}^{\mathrm{T}} = [\mathbf{q}_{b1,a+}^{\mathrm{T}} \quad \mathbf{q}_{b2,a+}^{\mathrm{T}} \quad \cdots \quad \mathbf{q}_{br,a+}^{\mathrm{T}}], \tag{121}$$

$\mathbf{q}_{b1,a+}^{\mathrm{T}}$, $\mathbf{q}_{b2,a+}^{\mathrm{T}}$, … and $\mathbf{q}_{br,a+}^{\mathrm{T}}$ are the vectors of Fourier coefficients corresponding to the elements of the first, second, … and $r$-th rows of $\bar{\mathbf{u}}(a, x_2)$, or we may call them the vectors of Fourier coefficients of displacement type boundary conditions, and $\mathbf{R}_{b,a+}$ is the matrix of Fourier coefficients corresponding to the matrix $\mathbf{\Gamma}(a, x_2)$.

Similarly, for the boundaries $x_1 = -a$, $x_2 = b$ and $x_2 = -b$, we can obtain three other sets of equations. Combining them with Eq. (120), we have

$$\mathbf{q}_b = \mathbf{R}_b \cdot \mathbf{q}, \tag{122}$$

where the matrix

$$\mathbf{R}_b = \begin{bmatrix} \mathbf{R}_{b,a+} \\ \mathbf{R}_{b,a-} \\ \mathbf{R}_{b,b+} \\ \mathbf{R}_{b,b-} \end{bmatrix}, \tag{123}$$

the vector of Fourier coefficients of displacement type boundary conditions
$$\mathbf{q}_b^{\mathrm{T}} = [\mathbf{q}_{b,a+}^{\mathrm{T}} \quad \mathbf{q}_{b,a-}^{\mathrm{T}} \quad \mathbf{q}_{b,b+}^{\mathrm{T}} \quad \mathbf{q}_{b,b-}^{\mathrm{T}}]. \tag{124}$$

We write the vectors of undetermined constants
$$\mathbf{q}_{03}^{\mathrm{T}} = [\mathbf{q}_0^{\mathrm{T}} \quad \mathbf{q}_3^{\mathrm{T}}], \tag{125}$$
$$\mathbf{q}_{12}^{\mathrm{T}} = [\mathbf{q}_1^{\mathrm{T}} \quad \mathbf{q}_2^{\mathrm{T}}], \tag{126}$$
and write the matrix $\mathbf{R}_b$ in the block matrix form
$$\mathbf{R}_b = [\mathbf{R}_{b,03} \quad \mathbf{R}_{b,12}], \tag{127}$$
then we can rewrite Eq. (122) as
$$\mathbf{q}_b = \mathbf{R}_{b,03} \cdot \mathbf{q}_{03} + \mathbf{R}_{b,12} \cdot \mathbf{q}_{12}, \tag{128}$$
and accordingly obtain
$$\mathbf{q}_{12} = -\mathbf{R}_{b,12}^{-1}\mathbf{R}_{b,03} \cdot \mathbf{q}_{03} + \mathbf{R}_{b,12}^{-1} \cdot \mathbf{q}_b. \tag{129}$$

Substitute Eq. (129) into Eq. (112) in [81], then we have a new representation of the Fourier series multiscale solution as
$$u = \mathbf{\Phi}_R^{\mathrm{T}} \cdot \mathbf{q}_R, \tag{130}$$
where the transformed matrix of basis functions
$$\mathbf{\Phi}_R^{\mathrm{T}} = \begin{bmatrix} \mathbf{\Phi}_{03}^{\mathrm{T}} & \mathbf{\Phi}_{12}^{\mathrm{T}} \end{bmatrix} \begin{bmatrix} \mathbf{I} & \mathbf{0} \\ -\mathbf{R}_{b,12}^{-1}\mathbf{R}_{b,03} & \mathbf{R}_{b,12}^{-1} \end{bmatrix}, \tag{131}$$
with the vectors of basis functions
$$\mathbf{\Phi}_{03}^{\mathrm{T}} = \begin{bmatrix} \mathbf{\Phi}_0^{\mathrm{T}} & \mathbf{\Phi}_3^{\mathrm{T}} \end{bmatrix}, \tag{132}$$
$$\mathbf{\Phi}_{12}^{\mathrm{T}} = \begin{bmatrix} \mathbf{\Phi}_1^{\mathrm{T}} & \mathbf{\Phi}_2^{\mathrm{T}} \end{bmatrix}, \tag{133}$$

and the transformed vectors of undetermined constants
$$\mathbf{q}_R^{\mathrm{T}} = [\mathbf{q}_{03}^{\mathrm{T}} \quad \mathbf{q}_b^{\mathrm{T}}]. \tag{134}$$

## 6. Derivation of the supplementary solution

For the linear differential equation with constant coefficients, the characteristic of the external force function $f$ has effects to some extent on the convergence and computational accuracy of the Fourier series multiscale solution. Under specific conditions, we can decompose the external force function into two parts, namely, the coarse scale component $f_s$ and the fine scale component $f - f_s$, to which the responses of the differential equation are determined with different methods.

For example, consider the two-dimensional full-range expansion in the domain $[-a,a] \times [-b,b]$, and we select a function $\varphi_s$ of a specific form (e.g. the algebraical polynomial) such that
$$\mathcal{L}\varphi_s = f_s. \tag{135}$$

Substitute it into Eq. (24) and then we have
$$\mathcal{L}(u - \varphi_s) = f - f_s, \quad (x_1, x_2) \in (-a,a) \times (-b,b). \tag{136}$$

Again, consider the boundary $x_1 = a$, and the corresponding boundary conditions are given by
$$\mathbf{B}(u - \varphi_s) = \mathbf{g} - \mathbf{B}\varphi_s. \tag{137}$$

Without loss of generality, suppose that the Fourier series multiscale solution of Eq.

(136), together with the prescribed boundary conditions in Eq. (137), is

$$u - \varphi_s = \varphi_0 + \varphi_1 + \varphi_2 + \varphi_3, \tag{138}$$

accordingly, the solution of Eq. (24) is

$$u = \varphi_0 + \varphi_1 + \varphi_2 + \varphi_3 + \varphi_s. \tag{139}$$

For brevity, we also call Eq. (139) the Fourier series multiscale solution of Eq. (24). Hence, the Fourier series multiscale solution includes not only the necessary terms, $\varphi_1$, $\varphi_2$ and $\varphi_0 + \varphi_3$, but also the optional and supplementary term $\varphi_s$. With the given definitions of the general solution and the particular solution of Eq. (24), we can further denote the function $\varphi_s$ by the supplementary solution of Eq. (24).

It is observed that reasonable choice of the supplementary solution provides an effective approach for improvement of convergence and computational accuracy of the Fourier series multiscale solution.

## 7. Derivation of discrete equations

After obtaining the representation of the Fourier series multiscale solution, we turn to some specific techniques for deriving discrete equations of undetermined constants concerned in the series solution. These derivation techniques can be classified into three types. The first is the Fourier coefficient comparison method, which is special in the Fourier series method for differential equations. The second is the weighted residual methods (for instance the collocation method), which are in the category of the numerical methods in strong form. The third is the variational methods (for instance the minimum potential energy method), which are in the category of the numerical methods in weak form. The applications of these three types of techniques will be illustrated with some typical classes of differential equations. However, as shown in Table 1, much attention should be paid to the matching of the derivation technique for discrete equations and the representation of the Fourier series multiscale solution in a specific scheme for solving the differential equations.

Table 1: Scheme for solving the differential equation with respect to the two-dimensional full-range expansion.

| Derivation technique for discrete equations | Undetermined constants | Representation of Fourier series multiscale solution |
|---|---|---|
| Fourier coefficient comparison method | $\mathbf{q}_0$, $\mathbf{q}_1$, $\mathbf{q}_2$, $\mathbf{q}_3$ | Eq. (139) |
| Collocation method | $\mathbf{q}_0$, $\mathbf{q}_1$, $\mathbf{q}_2$, $\mathbf{q}_3$ | Eq. (139) |
| Minimum potential energy method | $\mathbf{q}_0$, $\mathbf{q}_b$, $\mathbf{q}_3$ | Eq. (130) |

## 8. Multiscale characteristic of the solution

As to the Fourier series multiscale solution, the structural decomposition of the solutions is virtually also a scale decomposition of the solution space. Herein we take the full-range expansion in the two-dimensional domain $[-a,a]\times[-b,b]$ as an example. Then, in the structural decomposition of the solution, the internal function $\varphi_0(x_1,x_2)$ is a two-dimensional full-range Fourier series in the domain $[-a,a]\times[-b,b]$ and corresponds to the primary scale (or middle scale) of the decomposition scales; the boundary function $\varphi_1(x_1,x_2)$ is a one-dimensional full-range Fourier series on the interval $[-b,b]$ along the $x_2$ direction, of which the Fourier coefficients are the interpolation functions of the homogeneous solutions of the ordinary differential equations of the variable $x_1$, and corresponds to the primary scale and the adaptive scale (or hidden scale related to the original differential equation) along the $x_1$ direction; the boundary function $\varphi_2(x_1,x_2)$ is a one-dimensional full-range Fourier series on the interval $[-a,a]$ along the $x_1$ direction, of which the Fourier coefficients are the interpolation functions of the homogeneous solutions of the ordinary differential equations of the variable $x_2$, and corresponds to the primary scale and the adaptive scale (or hidden scale related to the original differential equation) along the $x_2$ direction; The corner function $\varphi_3(x_1,x_2)$ is a two-dimensional algebraical polynomial in the domain $[-a,a]\times[-b,b]$ and corresponds to the global scale (or large scale) of the decomposition scales. The global, large scale, the primary, middle scale and the adaptive, hidden scale compose the decomposition structure of scales in the process of solving the differential equation. Accordingly, the particular solution subspace (constructed by the sum of internal function $\varphi_0(x_1,x_2)$ and corner function $\varphi_3(x_1,x_2)$) on the large scale and middle scale and the general solution space (constructed by the sum of boundary functions $\varphi_1(x_1,x_2)$ and $\varphi_2(x_1,x_2)$) on the hidden scale compose the multiscale decomposition structure of the solution space of the differential equation.

The Fourier series multiscale solution is capable of adjusting adaptively the scales of analysis. It enables the Fourier series multiscale method, namely, the composite Fourier series method where homogeneous solutions of the differential equation are adopted as interpolation functions, to be a practical, feasible, effective and accurate multiscale method for computational analysis in science and engineering.

## 9. Conclusions

Effectiveness and high accuracy are the prerequisites for application of the multiscale methods. In this paper, we perform thorough investigation on the analytical analysis of multiscale phenomena inherent in the 2*r*-th order linear differential equations with constant coefficients and subjected to general boundary conditions. It is concluded that

1. We develop a new formulation of the composite Fourier series method where homogeneous solutions of the differential equations are adopted as interpolation functions, and establish the theoretical framework of the Fourier series multiscale method.

2. We specify the decomposition structures of solutions of the differential equations, which results in perfect integration of consistency and flexibility of the Fourier series solution.

3. We detail the implementation schemes for application of the Fourier series multiscale method.

4. We illustrate the capability of adaptive adjustment of the analysis scales in the Fourier series multiscale solutions, and discover the internal mechanism of conciseness and high efficiency of the Fourier series multiscale method.

The research achievements above lead to the developments of both the researches on the multiscale methods and the Fourier series methods for differential equations.

## References

1. Chaudhuri RA. On the roles of complementary and admissible boundary constraints in Fourier solutions to the boundary value problems of completely coupled $r$th order PDEs. *Journal of Sound and Vibration* 2002; **251**: 261-313.
2. Timoshenko S, Woinowsky-Krieger S. *Theory of Plates and Shells* (Second Edition). McGraw-Hill Book Company: Singapore, 1959.
3. Zhang FF. *Theory of Elastic Thin Plates* (Second Edition). Science Press: Beijing, 1984 (in Chinese).
4. Gorman DJ. *Free Vibration Analysis of Rectangular Plates*. Elsevier: North Holland, 1982.
5. Gorman DJ. A comprehensive study of the free vibration of rectangular plates resting on symmetrically distributed uniform elastic edge supports. *Journal of Applied Mechanics* 1989; **56**: 893-899.
6. Gorman DJ. Accurate free vibration analysis of point supported Mindlin plates by the superposition method. *Journal of Sound and Vibration* 1999; **219**: 265-277.
7. Gorman DJ. Free vibration and buckling of in-plane loaded plates with rotational elastic edge support. *Journal of Sound and Vibration* 2000; **229**: 755-773.
8. Gorman DJ. Free vibration analysis of completely free rectangular plates by the Superposition-Galerkin method. *Journal of Sound and Vibration* 2000; **237**: 901-914.
9. Gorman DJ. Accurate analytical type solutions for the free in-plane vibration of clamped and simply supported rectangular plates. *Journal of Sound and Vibration* 2004; **276**: 311-333.
10. Gorman DJ. Highly accurate free vibration eigenvalues for the completely free orthotropic plate. *Journal of Sound and Vibration* 2005; **280**: 1095-1115.
11. Gorman DJ. Exact solutions for the free in-plane vibration of rectangular plates with two opposite edges simply supported. *Journal of Sound and Vibration* 2006; **294**: 131-161.
12. Gorman DJ. On use of the Dirac delta function in the vibration analysis of elastic structures. *International Journal of Solids and Structures* 2008; **45**: 4605-4614.
13. Gorman DJ. Accurate in-plane free vibration analysis of rectangular orthotropic plates. *Journal of Sound and Vibration* 2009; **323**: 426-443.
14. Gorman DJ, Ding W. Accurate free vibration analysis of laminated symmetric cross-ply rectangular plates by the superposition-Galerkin method. *Composite Structures* 1995; **31**: 129-136.
15. Gorman DJ, Ding W. Accurate free vibration analysis of clamped antisymmetric angle-ply laminated rectangular plates by the Superposition-Galerkin method. *Composite Structures* 1996; **34**: 387-395.
16. Gorman DJ, Ding W. Accurate free vibration analysis of the completely free rectangular Mindlin plate. *Journal of Sound and Vibration* 1996; **189**: 341-353.
17. Gorman DJ, Ding W. the Superposition-Galerkin method for free vibration analysis of rectangular plates. *Journal of Sound and Vibration* 1996; **194**: 187-198.
18. Gorman DJ, Ding W. Accurate free vibration analysis of completely free symmetric cross-ply rectangular laminated plates. *Composite Structures* 2003; **60**: 359-365.
19. Gorman DJ, Garibaldi L. Accurate analytical type solutions for free vibration frequencies and mode shapes of multi-span bridge decks: the span-by-span approach. *Journal of Sound and Vibration* 2006;

**290**: 321-336.

20. Gorman DJ, Singhal R. Free vibration analysis of cantilever plates with step discontinuities in properties by the method of superposition. *Journal of Sound and Vibration* 2002; **253**: 631-652.

21. Gorman DJ, Singhal R. Steady-state response of a cantilever plate subjected to harmonic displacement excitation at the base. *Journal of Sound and Vibration* 2009; **323**: 1003-1015.

22. Gorman DJ, Yu SD. A review of the superposition method for computing free vibration eigenvalues of elastic structures. *Computers and Structures* 2012; **104-105**: 27-37.

23. Huang Y. *Theory of Thin Plates*. Press of National University of Defense Technology: Changsha, 1992 (in Chinese).

24. Huang Y, Tan YZ, Xu XL. A general analytical solution of free vibration of orthotropic rectangular thin plates. *Engineering Mechanics* 2001; **18**: 45-52 (in Chinese).

25. Huang Y, Jiang YQ, Li JW. Static analysis of anisotropic plates on elastic foundation. *Chinese Journal of Applied Mechanics* 2006; **23**: 466-469 (in Chinese).

26. Huang Y, Lei YJ, Shen HJ. Free vibration of anisotropic rectangular plates by general analytical method. *Applied Mathematics and Mechanics* (English Edition) 2006; **27**: 461-467.

27. Huang Y, Liao Y, Xie Y. Free vibration of compressed orthotropic plate on two parameter elastic foundation. *Engineering Mechanics* 2006; **23**: 46-49 (in Chinese).

28. Huang Y, Yuan DC. Static analysis of symmetric angle-ply laminated plates by analytical method. *American Institute of Aeronautics and Astronautics Journal* 2006; **44**: 667-669.

29. Yang DS, Pan J, Huang Y. A general solution of an isotropic thin plate in bending problem. *Chinese Journal of Computational Mechanics* 2002; **19**: 286-290 (in Chinese).

30. Yang DS, Huang Y, Ren XH. Analysis of symmetric laminated rectangular plates in plane stress. *Applied Mathematics and Mechanics* (English Edition) 2006; **27**: 1719-1726.

31. Yang DS, Huang Y, Li GL. Shear buckling analysis of anisotropic rectangular plates. *Chinese Journal of Applied Mechanics* 2012; **29**: 220-224 (in Chinese).

32. Zheng RY, Huang Y, Liao YH. Free vibration analysis of rectangular plates with mixed boundary conditions. *Engineering Mechanics* 2008; **25**: 13-17 (in Chinese).

33. Zheng RY, Huang Y, Li GL. Vibration analysis of rectangular plates with general elastic boundary supports. *Journal of Mechanical Strength* 2010; **32**: 324-328 (in Chinese).

34. Zheng RY, Lin WF, Huang Y. Stability analysis of anisotropic rectangular plates. *Engineering Mechanics* 2009; **26**: 30-33 (in Chinese).

35. Chaudhuri RA, Abu-Arja KR. Exact solution of shear-flexible doubly curved anti-symmetric angle-ply shells. *International Journal of Engineering Science* 1988; **26**: 587-604.

36. Chaudhuri RA, Abu-Arja KR. Static analysis of moderately thick anti-symmetric angle-ply cylindrical panels and shells. *International Journal of Solids and Structures* 1991; **28**: 1-16.

37. Chaudhuri RA, Balaraman K, Kunukkasseril VX. A combined theoretical and experimental investigation on free vibration of thin symmetrically laminated anisotropic plates. *Composite Structures* 2005; **67**: 85-97.

38. Chaudhuri RA, Kabir HRH. A boundary discontinuous Fourier solution for clamped transversely isotropic (pyrolytic graphite) Mindlin plates. *International Journal of Solids and Structures* 1992; **30**: 287-297.

39. Chaudhuri RA, Kabir HRH. Fourier solution to higher-order theory based laminated shell boundary-value problem. *American Institute of Aeronautics and Astronautics Journal* 1995; **33**: 1681-1688.

40. Chaudhuri RA, Kabir HRH. Effect of boundary constraint on the frequency response of moderately thick doubly curved cross-ply panels using mixed fourier solution functions. *Journal of Sound and Vibration* 2005; **283**: 263-293.
41. Kabir HRH. A novel approach to the solution of shear flexible rectangular plates with arbitrary laminations. *Composite: Part B* 1996; **27B**: 95-104.
42. Kabir HRH. On boundary value problems of moderately thick shallow cylindrical panels with arbitrary laminations. *Composite Structures* 1996; **34**: 169-184.
43. Kabir HRH. Anti-symmetric angle-ply laminated thick cylindrical panels. *International Journal of Solids and Structures* 1998; **35**: 3717-3735.
44. Kabir HRH. Free vibration response of shear-deformable antisymmetric cross-ply cylindrical panels. *Journal of Sound and Vibration* 1998; **217**: 601-618.
45. Kabir HRH. Application of linear shallow shell theory of Reissner to frequency response of thin cylindrical panels with arbitrary lamination. *Composite Structures* 2002; **56**: 35-52.
46. Kabir HRH. On free vibration response and mode shapes of arbitrarily laminated rectangular plates. *Composite Structures* 2004; **65**: 13-27.
47. Kabir HRH, Al-Khaleefi AM, Al-Marzouk M. Double orthogonal set of solution functions for cross-ply laminated shear flexible cylindrical/doubly curved panels. *Composite Structures* 2003; **59**: 189-198.
48. Kabir HRH, Al-Khaleefi AM, Chaudhuri RA. Free vibration analysis of thin arbitrarily laminated anisotropic plates using boundary-continuous displacement Fourier approach. *Composite Structures* 2001; **53**: 469-476.
49. Kabir HRH, Al-Khaleefi AM, Chaudhuri RA. Frequency response of a moderately thick antisymmetric cross-ply cyclindrical panel using mixed type of fourier solution functions. *Journal of Sound and Vibration* 2003; **259**: 809-828.
50. Kabir HRH, Chaudhuri RA. Boundary-continuous Fourier solution for clamped Mindlin plates. *Journal of Engineering Mechanics* 1992; **118**: 1457-1467.
51. Kabir HRH, Hamad MAM, Al-Duaij J, Mariam JJ. Thermal buckling response of all-edge clamped rectangular plates with symmetric angle-ply lamination. *Composite Structures* 2007; **79**: 148-155.
52. Oktem AS. The effect of boundary conditions on the response of laminated thick composite plates and shells. *Ph.D. Thesis*, The University of Utah, 2005.
53. Oktem AS, Alankaya V, Soares C.G. Boundary-discontinuous Fourier analysis of simply supported cross-ply plates. *Applied Mathematical Modelling* 2013; **37**: 1378-1389.
54. Oktem AS, ChaudhuriRA. Levy type analysis of cross-ply plates based on higher-order theory. *Composite Structures* 2007; **78**: 243-253.
55. Oktem AS, Chaudhuri RA. Fourier solution to a thick cross-ply Levy type clamped plate problem. *Composite Structures* 2007; **79**: 481-492.
56. Oktem AS, Chaudhuri RA. Levy type Fourier analysis of thick cross-ply doubly curved panels. *Composite Structures* 2007; **80**: 475-488.
57. Oktem AS, Chaudhuri RA. Fourier analysis of thick cross-ply Levy type clamped doubly-curved panels. *Composite Structures* 2007; **80**: 489-503.
58. Oktem AS, Chaudhuri RA. Boundary discontinuous Fourier analysis of thick cross-ply clamped plates. *Composite Structures* 2008; **82**: 539-548.
59. Oktem AS, Chaudhuri RA. Effect of inplane boundary constraints on the response of thick general (unsymmetric) cross-ply plates. *Composite Structures* 2008; **83**: 1-12.

60. Oktem AS, Chaudhuri RA. Higher-order theory based boundary-discontinuous Fourier analysis of simply supported thick cross-ply doubly curved panels. *Composite Structures* 2009; **89**: 448-458.
61. Oktem AS, Soares CG. Boundary discontinuous Fourier solution for plates and doubly curved panels using a higher order theory. *Composites: Part B* 2011; **42**: 842-850.
62. Yan ZD. *The Fourier Series Method in Structural Mechanics*. Tianjin University Press: Tianjin, 1989 (in Chinese)
63. Li WL. Free vibrations of beams with general boundary conditions. *Journal of Sound and Vibration* 2000; **237**: 709-725.
64. Li WL. Dynamic analysis of beams with arbitrary elastic supports at both ends. *Journal of Sound and Vibration* 2001; **246**: 751-756.
65. Li WL. Comparision of Fourier sine and cosine series expansions for beams with arbitrary boundary conditions. *Journal of Sound and Vibration* 2002; **255**: 185-194.
66. Li WL, Daniels M. A Fourier series method for the vibrations of elastically restrained plates arbitrarily loaded with springs and masses. *Journal of Sound and Vibration* 2002; **252**: 768-781.
67. Li WL. Vibration analysis of rectangular plates with general elastic boundary supports. *Journal of Sound and Vibration* 2004; **273**: 619-635.
68. Li WL, Xu HA. An exact Fourier series method for the vibration analysis of multispan beam systems. *Journal of Computational and Nonlinear Dynamics* 2009; **4**: 021001.
69. Li WL, Zhang XF, Du JT, Liu ZG. An exact series solution for the transverse vibration of rectangular plates with general elastic boundary supports. *Journal of Sound and Vibration* 2009; **321**: 254-269.
70. Li WL, Bonilha MW. Vibrations of two beams elastically coupled together at an arbitrary angle. *Acta Mechanica Solida Sinica* 2012; **25**: 61-72.
71. Khov H, Li WL, Gibson RF. An accurate solution method for the static and dynamic deflections of orthotropic plates with general boundary conditions. *Composite Structures* 2009; **90**: 474-481.
72. Du JT, Li WL, Liu ZG, Yang TJ, Jin GY. Free vibration of two elastically coupled rectangular plates with uniform elastic boundary restraints. *Journal of Sound and Vibration* 2011; **330**: 788-804.
73. Dai L, Yang TJ, Du JT, Li WL, Brennan MJ. An exact series solution for the vibration analysis of cylindrical shells with arbitrary boundary conditions. *Applied Acoustics* 2013; **74**: 440-449.
74. Zhang XF, Li WL. Vibrations of rectangular plates with arbitrary non-uniform elastic edge restraints. *Journal of Sound and Vibration* 2009; **326**: 221-234.
75. Zhang XF, Li W.L. A unified approach for predicting sound radiation from baffled rectangular plates with arbitrary boundary conditions. *Journal of Sound and Vibration* 2010; **329**: 5307-5320.
76. Xu HG, Du JT, Li WL. Vibrations of rectangular plates reinforced by any number of beams of arvitrary lengths and placement angles. *Journal of Sound and Vibration* 2010; **329**: 3759-3779.
77. Onate E, Miquel J, Hauke G. Stabilized formulation for the advection-diffusion-absorption equation using finite calculus and linear finite elements. *Computer Methods in Applied Mechanics and Engineering* 2006; **195**: 3926-3946.
78. Parvazinia M, Nassehi V. Multiscale finite element modeling of diffusion-reaction equation using bubble functions with bilinear and triangular elements. *Computer Methods in Applied Mechanics and Engineering* 2007; **196**: 1095-1107.
79. Parvazinia M, Nassehi V, Wakeman RJ. Multi-scale finite element modelling of laminar steady flow through highly permeable porous media. *Chemical Engineering Science* 2006; **61**: 586-596.
80. Parvazinia M, Nassehi V, Wakeman RJ. Multi-scale finite element modelling using bubble function method for a convection-diffusion problem. *Chemical Engineering Science* 2006; **61**: 2742-2751.

81. Sun WM, Zhang ZM. Fourier series (based) multiscale method for computational analysis in science and engineering: II. Composite Fourier series method for simultaneous approximation to functions and their (partial) derivatives. doi:10.48550/ARXIV.2208.03749.